\documentclass[a4j,10pt]{article}
\usepackage{graphicx}
\usepackage{latexsym}
\usepackage{amsmath}
\usepackage{amsthm}
\usepackage{amssymb}
\usepackage{color}
\newtheorem{thm}{Theorem}
\newtheorem{prop}[thm]{Proposition}

\newtheorem{claim}{Claim}

\newtheorem{lm}[thm]{Lemma}

\newtheorem{conj}[thm]{Conjecture}

\newcounter{Case}
\newcounter{Subcase}[Case]
\def\komento#1{\ }
\def\qed{\hfill$\Box$}
\def\ed{d_H^e}

\def\cl{{\mbox{\textit{cl}}}}

\allowdisplaybreaks
\begin{document}
\hspace*{0.5cm}
\begin{center}
{\Large
A supplement to
``Induced nets and Hamiltonicity of claw-free graphs''}
\bigskip\\
{\large
Shuya Chiba\footnotemark[1] \footnotemark[2]\qquad
Jun Fujisawa\footnotemark[3] \footnotemark[4]\\
\footnotetext[1]{%
Applied Mathematics, Faculty of Advanced Science and Technology,
Kumamoto University,
2-39-1, Kurokami, Kumamoto 860--8555,
Japan.
\texttt{schiba@kumamoto-u.ac.jp}
}%
\footnotetext[2]{%
work supported by Japan Society for the Promotion of Science,
Grant-in-Aid for Scientific Research (C) 17K05347}%
\footnotetext[3]{%
Faculty of Business and Commerce,
Keio University,
Hiyoshi 4--1--1, Kohoku-Ku,
Yokohama, Kanagawa 223--8521,
Japan.
\texttt{fujisawa@fbc.keio.ac.jp}
}%
\footnotetext[4]{%
work supported by Japan Society for the Promotion of Science,
Grant-in-Aid for Scientific Research (B) 16H03952 and
Grant-in-Aid for Scientific Research (C) 17K05349}%
}
\end{center}
\medskip
\begin{abstract}
This note supplements our paper
``Induced nets and Hamiltonicity of claw-free graphs'',
by giving the detailed proof that were omitted in it.
\end{abstract}
\bigskip
\bigskip
%
%
%
%
%
%
In \cite{conj}, the following conjecture is proposed.
\begin{conj}[Broersma \cite{conj}]
Let $G$ be a $2$-connected claw-free graph of order $n$.  
If every endvertex of each induced net in $G$ has degree
at least $\frac{n-2}3$, 
then $G$ is hamiltonian.
\label{conj}
\end{conj}
In \cite{CF},
the authors announced that
Conjecture \ref{conj} is proved in the affirmative,
however we omitted the detailed proof for the graphs of order at most $32$.
The aim of this note is to give a complete proof of this case.

This note is a supplement to \cite{CF}.
For terminology and notation not defined here,
we refer the readers to \cite{CF} and \cite{Ds}.
Let us start with preliminary results.
%
%
%
%
%
%
\begin{thm}[Matthews and Sumner \cite{MS}]
Let $G$ be a $2$-connected claw-free graph of order $n$. 
If minimum degree of $G$ is at least $\frac{n-2}3$, then $G$ is hamiltonian.
\label{MS}
\end{thm}
\begin{thm}[Ryj\'a\v{c}ek \cite{Ry}]
\label{Ry}
Let $G$ be a claw-free graph.
Then $\cl(G)$ is uniquely defined and is the line graph of some triangle-free simple graph.
Moreover, $G$ is hamiltonian if and only if $\cl(G)$ is hamiltonian.
\end{thm}
\begin{thm}[Harary and Nash-Williams \cite{HN}]
Let $H$ be a multigraph with $|E(H)| \ge 3$.
Then the line graph $L(H)$ is hamiltonian
if and only if $H$ has a DCT.
\label{HN}
\end{thm}
\begin{prop}[Catlin \cite{Ca}]
Let $H$ be a multigraph and $F \subset H$ be a collapsible subgraph.
\begin{itemize}
\item[\normalfont{i)}]
If $H/F$ has a DCT containing $v_F$,
then $H$ has a DCT containing all the vertices in $F$.
\item[\normalfont{ii)}]
If $H/F$ is collapsible, then $H$ is collapsible.
\end{itemize}
\label{col}
\end{prop}
\begin{thm}[Lai \cite{Lai-cycle}]
Let $G$ be a $2$-connected graph with minimum degree at least $3$.
If every edge of $G$ lies in a cycle of length at most $4$,
then $G$ is collapsible.
\label{Lai-cycle}
\end{thm}
\begin{lm}[Chiba and Fujisawa \cite{CF}]
Let $G$ be a $2$-connected claw-free graph of order at least $3$
such that
every endvertex of each induced net in $G$ has degree
at least $\frac{|V(G)|-2}3$. 
Let $H$ be the triangle-free graph such that $L(H) = \cl(G)$,
and let $\Lambda$ be the subdivided claw of $H$ such that
$V(\Lambda) = \{R_0, R_1, R_2, R_3, R_1^+, R_2^+, R_3^+\}$
and $R_0 R_i R_i^+$ is an induced path of $\Lambda$
for $i = 1,2,3$.
Moreover,
let $R_0 R_i = x_i$ and $R_i R_i^+ = y_i$ for $i=1,2,3$.
Then the following holds.
\begin{itemize}
\item[\normalfont{i)}]
There exists an induced net $N^0 = N(x_1', x_2', x_3'; y_1',y_2',y_3')$ in $G$
such that $y_i' \in  \{y_i, x_i\} \cup l_H(R_i) \cup l_H(R_0)$.
\item[\normalfont{ii)}]
Each $y_i'$ is heavy. Moreover, if $y_i' \in \{x_i\} \cup l_H(R_0)$, then
$\ed(y_i') \ge \frac{|E(H)|-2}3 + 2 + |J|$,  
where $J = \{y_j' \mid j \neq i, \ y_j' \in \{x_j\} \cup l_H(R_0)\}$.
\item[\normalfont{iii)}]
$x_1 x_2 x_3$ is a triangle of $G$ if and only if 
$y_i' \in \{y_i\} \cup l_H(R_i)$ for each $i$.
\end{itemize}
\label{endmove}
\end{lm}
\begin{lm}[Chiba and Fujisawa \cite{CF}]
Let $H$ be an essentially $2$-connected multigraph
and let $\Xi$ be a collapsible subgraph of $H$.
If $|E(H-\Xi)| \le 3$,
then there exists a DCT of $H$.
\label{nokorisanhen_c}
\end{lm}
\begin{lm}[Chiba and Fujisawa \cite{CF}]
Both of $K_{3,3}$ and $K_{3,3}^-$ is collapsible.
\qed
\label{lm_k33-}%
\end{lm}%
\begin{lm}[Chiba and Fujisawa \cite{CF}]
Let $H$ be an essentially $2$-edge-connected triangle-free simple graph
which does not contain a DCT
and let $n = |E(H)|$.
If $\{e_1, e_2, e_3\}$ is a matching of $H$,
then $\sum_{i=1}^3 \ed(e_i) \le n+1$.
\label{jisuukeisan}
\end{lm}
\begin{thm}[Chiba and Fujisawa \cite{CF}]
Let $G$ be a graph of order $n$
which satisfies the assumption of Conjecture \ref{conj}
and let $H$ be the triangle-free graph such that
$L(H) = \cl(G)$.
Then there exists either a DCT
or a heavy matching of size $4$ in $H$.
\label{main}
\end{thm}
%
%
%
%
%
%
%
For a graph $G$ and $X,Y \subseteq V(G)$ with $X \cap Y = \emptyset$,
the set of edges between $X$ and $Y$ is denoted by $E_G(X,Y)$,
or simply $E(X,Y)$.
Moreover, we use $e(X,Y) = |E(X,Y)|$.
If there is no fear of confusions,
we often identify a vertex $v$ and $\{v\}$,
and a subgraph $F$ of $G$ and $V(H)$.
For example, we use $e(v,F)$ instead of $e(\{v\}, V(F))$.
\begin{lm}
Let $G$ be a triangle-free graph with $|V(G)|=7$
which contains a cycle $C$ of length $5$,
and let $V(G) \setminus V(C) = \{w_1, w_2\}$.
If $e(w_i, C) = 2$ for $i=1,2$,
then $G$ has a spanning closed trail
which contains all the edges of $E(\{w_1, w_2\},C)$.
In addition, if $w_1 w_2 \in E(G)$,
then $G$ is collapsible.
\label{5cycle-lm}
\end{lm}
\noindent
\textit{Proof.}\quad
Let  $C = u_1 u_2 \ldots u_5 u_1$.
Since $G$ is triangle-free,
$C$ has no chord, and by symmetry we may assume that
$N(w_1) \cap V(C) = \{u_1, u_3\}$.
Then we may assume
$N(w_2) \cap V(C) = \{u_1, u_3\}$, $\{u_2, u_4\}$ or $\{u_3, u_5\}$.
Let $E^- = \emptyset$ (resp. $\{u_1 u_2, u_3 u_4\}$ and $\{u_5 u_1\}$)
in the first (resp. second and last) case,
then $( E(C) \setminus E^- ) \cup E(\{w_1, w_2\},C) $ induces the required closed trail.

If $w_1 w_2 \in E(G)$, then we may assume that
$N(w_2) \cap V(C) = \{u_2, u_4\}$
since $G$ is triangle-free.
Let $G' = G - \{u_5\}$,
then $G' \simeq K_{3,3}^-$
and hence $G'$ is collapsible by Lemma \ref{lm_k33-}.
Moreover, since the edges $G/G'$ form multiple edges,
$G/G'$ is collapsible.
By Proposition \ref{col} ii),
$G$ is collapsible.
\qed

\begin{lm}
Let $G$ be an essentially $2$-edge-connected graph and
let $xy \in E(G)$ such that $d_G(x), d_G(y) \ge 2$.
If $|E(G-\{x,y\})| \le 2$,
then there exists a DCT of $G$
containing $x$ and $y$.
\label{nokori2hen}
\end{lm}
\noindent
\textit{Proof.}\quad
By assumption,
we can take a cycle $C$ which contains the edge $xy$ in $G$.
Let $E_0 = E(G-\{x,y\})$.
If $E_0 \subseteq E(C)$,
then $C$ is a DCT of $G$.
Thus we assume $E_0 \not \subseteq E(C)$.
Since $|E_0| \le 2$, we have $|V(C)| \le 4$.
If $|V(C)| = 3$,
then $C$ is collapsible,
and hence we obtain a desired closed trail
by Lemma \ref{nokorisanhen_c}.
Hence we assume $|V(C)| = 4$.
Then the edge, say $e_1$, of $C-\{x, y\}$ is in $E_0$.
Let $e_2$ be the other edge of $E_0$,
then we may assume that $e_2 \in E(G-V(C))$,
since otherwise $C$ is a DCT of $G$.
Note that $E(G-\{x,y\}) = \{e_1, e_2\}$,
which implies
$N(v) \setminus V(e_{i}) \subseteq \{x,y\}$ 
for each 
$v \in V(e_i)$ $(i = 1, 2)$.
Since $G$ is essentially $2$-edge-connected,
there exist two edges between $V(e_2)$ and $\{x, y\}$,
and thus there exists a cycle $C'$ in $G[V(e_2) \cup \{x,y\}]$.
Then $C \triangle C'$ is a desired closed trail,
where $\triangle$ denotes the symmetric difference.
\qed
\medskip

We can now prove the main result.
\begin{thm}
Let $G$ be a $2$-connected claw-free graph of order $n \le 32$.
If every endvertex of each induced net in $G$ has degree
at least $\frac{n-2}3$, then $G$ is hamiltonian.
\label{thn}
\end{thm}
\noindent\textit{Proof.}\quad
Let $H$ be the triangle-free graph such that
$L(H) = \cl(G)$.
By Theorems \ref{Ry} and \ref{HN},
it suffices to prove that
$H$ has a DCT.
\medskip

\noindent
\textbf{Case 1.}\quad
$n \ge 15$.
\medskip

By Theorem \ref{main},
we may assume that there exists a heavy matching $M$ of size $4$ in $H$.
Let
$\Xi^* = H[V(M)]$ and
$E_0 = E(H - V(\Xi^*))$.
Since $\Xi^*$ is triangle-free,
it follows from Tur{\'a}n's theorem that $|E(\Xi^*)| \le 16$.
Moreover, since
$n  =  \sum_{e \in M} \ed(e) + |M| - |E(\Xi^*) \setminus M| + |E_0|
 \ge  4 \cdot \frac{n-2}3 + 4 - (|E(\Xi^*)|-4) +|E_0|$,
we have
\begin{eqnarray}
|E_0| \le |E(\Xi^*)| - \frac{n+16}3.
\label{Dec14}
\end{eqnarray}
Since $n \ge 15$ and $|E(\Xi^*)| \le 16$,
(\ref{Dec14}) yields
\begin{eqnarray}
|E_0| \le |E(\Xi^*)| - 11
\label{Dec8_1}
\end{eqnarray}
and
\begin{eqnarray}
|E_0| \le \frac{32 - n}3.
\label{Dec8_2}
\end{eqnarray}

Since $H$ is essentially $2$-edge-connected,
we have $|E(H) \setminus \left( E(\Xi^*) \cup E_0\right) | \ge 2$
if $E_0 \neq \emptyset$
(consider two edge-disjoint paths joining an edge of $E_0$ and $\Xi^*$).
This implies
\begin{eqnarray}
|E_0| \le \max\{0, n - |E(\Xi^*)| - 2\}.
\label{Sep29_2}
\end{eqnarray}

Assume, 
for the moment that,
$|E(\Xi^*)| \ge 14$.
Let $E(\Xi^*) = \{e_1, e_2, e_3, e_4\}$ and
$f_{ij} = e(V(e_i), V(e_j))$
for $1 \le i < j \le 4$.
Then $\sum_{1 \le i < j \le 4} f_{ij} \ge |E(\Xi^*)-M| \ge 10$.
By the pegionhole principle,
we can take
$i, j, k$ with $1 \le i < j < k \le 4$
so that
$f_{ij} + f_{ik} + f_{jk} \ge 5$. 
Let $\Xi = H[V(e_i) \cup V(e_j) \cup V(e_k)]$,
then 
$|E(\Xi)| = 3 + f_{ij} + f_{ik} + f_{jk} \ge 8$.
Since $|V(\Xi)| = 6$ and $\Xi$ is triangle-free,
we can deduce that $\Xi$ is bipartite.
Hence $\Xi \simeq K_{3,3}$ or $K_{3,3}^-$,
and by Theorem \ref{Lai-cycle} and Lemma \ref{lm_k33-},
$\Xi$ is collapsible.
Furthermore, since $|E(\Xi^*)| \ge 14$ and $|E(\Xi)| \le 9$,
we can deduce that
both of the endvertices of the edge of $\Xi^* - V(\Xi)$ have a neighbor in $\Xi$.
Hence $\Xi^*/\Xi$ is a triangle with some multiple edges.
By Proposition \ref{col} ii), $\Xi^*$ is collapsible.

Recall that $v_{\Xi^*}$ is the vertex in $H/\Xi^*$ to which $\Xi^*$ contracts,
and $|E(H/\Xi^* - \{v_{\Xi^*}\})| = |E_0|$.
If $|E_0| \le 3$, then by Lemma \ref{nokorisanhen_c}
there exists a DCT of $H$.
Hence we may assume that $|E_0| \ge 4$.
Then (\ref{Dec8_1}) yields $|E(\Xi^*)| \ge 15$ and
(\ref{Dec8_2}) yields $n \le 20$.
This contradicts (\ref{Sep29_2}),
and thus we have $|E(\Xi^*)| \le 13$.
By (\ref{Dec8_1}), we have $|E_0| \le 2$.
We distinguish two cases.
\medskip\\
\textbf{Case 1.1.}\quad
$\Xi^*$ is bipartite.
\medskip

First assume that $12 \le |E(\Xi^*)| \le 13$.
If $\Xi^*$ has no vertex of degree one,
then we can deduce that either $\Xi^*$ has minimum degree at least $3$
or $\Xi^*$ contains a subgraph $\Xi'$ such that $\Xi' \simeq K_{3,3}$ or $K_{3,3}^-$.
In the former case,
$\Xi^*$ satisfies the assumption of Theorem \ref{Lai-cycle},
and in the latter case, $\Xi'$ is collapsible and
$\Xi^*/\Xi'$ is a triangle (possibly with some multiple edges)
or a connected graph of order $3$ consisting of multiple edges.
Thus $\Xi^*$ is collapsible in either case.
By Lemma \ref{nokorisanhen_c},
$H$ has a DCT, since $|E_0| \le 2$.
Therefore we may assume that
$\Xi^*$ has a vertex, say $x$, of degree one.
Let $\Xi' = \Xi^*-\{x\}$.
Then $|E(\Xi')| = |E(\Xi^*)|-1  \ge 11$.
Since $\Xi' \subseteq K_{3,4}$,
$\Xi'$ contains a subgraph $\Xi''$ which is isomorphic to $K_{3.3}$.
By Theorem \ref{Lai-cycle}, $\Xi''$ is collapsible,
and since $\Xi'/\Xi''$ consists of multiple edges with two or three edges,
it follows from Proposition \ref{col} ii) that $\Xi'$ is collapsible.
Note that $\Xi^*/\Xi' \simeq K_2$ with $V(\Xi^*/\Xi') = \{x, v_{\Xi'}\}$.
If $x$ has degree one in $H$,
then by Lemma \ref{nokorisanhen_c},
$H$ has a DCT,
since $|E(H-V(\Xi'))| = |E_0| \le 2$. 
Hence we assume that $x$ has degree at least two in $H$.
Then, since $H$ is essentially $2$-edge-connected,
$v_{\Xi'}$ has degree at least two in $H/\Xi'$.
By Lemma \ref{nokori2hen},
there exists a DCT of $H/\Xi'$
which contains $v_{\Xi'}$ and $x$,
and hence by Proposition \ref{col} i),
$H$ has a DCT.

Next assume that $|E(\Xi^*)| \le 11$.
Then by (\ref{Dec8_1}), we have $|E(\Xi^*)| = 11$ and $E_0 = \emptyset$.
If the minimum degree of $\Xi^*$ is two,
then since $|E(\Xi^*)| = 11$,
we can deduce that $\Xi^*$ has a hamiltonian cycle,
which is a DCT of $H$.
Hence we may assume that $\Xi^*$ has a vertex $x$ of degree one.
Then $\Xi^* - \{x\}$ contains a subgraph $\Xi'$ which is isomorphic to
$K_{3.3}$ or $K_{3,3}^-$.
By Theorem \ref{Lai-cycle} and Lemma \ref{lm_k33-},
$\Xi'$ is collapsible.
Moreover, $\Xi^* / \Xi'$ is a path, say $P$, of length two (possibly with some multiple edges). 
Since $H$ is essentially $2$-edge-connected,
we can take a closed trail in $H/\Xi'$
which contains $V(P) \setminus V_1(H)$.
This is a DCT of $H/\Xi'$,
and hence by Proposition \ref{col} i),
$H$ has a DCT.
\medskip\\
\textbf{Case 1.2.}\quad
$\Xi^*$ is non-bipartite.
\medskip

If $\Xi^*$ has an induced cycle $C_0$ of length $7$,
then since $\Xi^*$ is triangle-free,
$e(x, C_0) \le 3$,
where $x$ is the vertex of $\Xi^*-V(C_0)$.
Thus $|E(\Xi^*)| \le 10$,
which contradicts (\ref{Dec8_1}).
Hence every cycle of length $7$ must have a chord.
Therefore, $\Xi^*$ has a cycle of length $5$,
say $C$.
Let $C=u_1 u_2 \ldots u_5 u_1$ and $\Xi^* - V(C) = W$.
Since $\Xi^*$ is triangle-free,
$C$ has no chord
and $|E(W)| \le 2$.
Moreover, $e(w, C) \le 2$ for any $w \in V(W)$.

First consider the case $|E(\Xi^*)| \ge 12$.
Recall that $|E_0| \le 2$.
By Lemma \ref{nokorisanhen_c}, it suffices to show that
$\Xi^*$ is collapsible.
If there exists an edge $w_1w_2$ in $W$
such that $e(w_i,C) =2$ for $i=1,2$,
then by Lemma \ref{5cycle-lm}, $\Xi^*[V(C) \cup \{w_1, w_2\}]$ is collapsible,
and since the vertex of $\Xi^*-(V(C) \cup \{w_1, w_2\})$ has degreeat least $2$,
$\Xi^*$ is collapsible.
Hence for each edge $w_1 w_2$ of $W$, we may assume that
$e(w_i,C) \le 1$ for $i=1$ or $2$.
Then,
since $|E(\Xi^*) \setminus E(C)|\ge 7$, we can deduce that
$\Xi^*[W]$ is a path, say $w_1w_2w_3$, of length two,
$e(w_i, C) = 2$ for $i=1,3$ and $e(w_2, C) = 1$.
Since $\Xi^*$ is triangle-free,
we may assume that
$N_{\Xi^*}(w_1) \cap V(C) = \{u_2, u_4\}$.
If $N(w_2) \cap V(C) = \{u_3\}$,
then by symmetry we may assume that
$N_{\Xi^*}(w_3) \cap V(C) = \{u_2, u_4\}$ or $\{u_4, u_1\}$.
In the former case,
let $\Xi' = \Xi^* - \{u_1, u_5\}$,
then $\Xi' \simeq K_{3,3}$ and thus it is collapsible.
Moreover, $\Xi^* / \Xi'$ is a triangle,
and hence $\Xi^*$ is collapsible by Proposition \ref{col} ii).
In the latter case,
let $C'$ be the cycle $u_1 u_2 u_3 w_2 w_3 u_1$
and $\Xi' = \Xi^* - \{u_5\}$.
Then $C'$ is a $5$-cycle, $w_1 u_4 \in E(\Xi')$
and each of $w_1$ and $u_4$ has two neighbors in $C'$.
Hence by Lemma \ref{5cycle-lm}, $\Xi'$ is collapsible,
and since $u_5$ has two neighbors in $\Xi'$,
$\Xi^*$ is collapsible.
If $N(w_2) \cap V(C) = \{u_1\}$,
then
by the symmetry of $w_{1}$ and $w_{3}$, we may assume that
$N_{\Xi^*}(w_3) \cap V(C) = \{u_2, u_4\}$ or $\{u_3, u_5\}$.
In the former case, $\Xi^* - \{u_3, u_5\}$ is collapsible
since $\Xi^* - \{u_3, u_5\} \simeq K_{3,3}^-$,
and hence $\Xi^*$ is collapsible
because both of $u_3, u_5$ have degree at least two in $\Xi^*$.
In the latter case, $\Xi^*$ satisfies the assumption of Theorem \ref{Lai-cycle},
and hence $\Xi^*$ is collapsible.

Next consider the case $|E(\Xi^*)| \le 11$.
Then by (\ref{Dec8_1}), we have $|E(\Xi^*)|=11$ and $E_0 = \emptyset$.
Hence it suffices to show that $\Xi^*$ has a spanning closed trail.
Let $V(W) = \{w_1, w_2, w_3\}$.
Without loss of generality, we may assume that
$E(W) = \{w_1 w_2\}$ if $|E(W)| =1$, and
$E(W) = \{w_1 w_2, w_2 w_3 \}$ if $|E(W)| =2$.
Let $\omega_i = e(w_i, C)$ for $1 \le i \le 3$.
By symmetry, we may assume that one of the following holds:
\begin{quote}
1)\ \ 
$|E(W)| =2$, $(\omega_1, \omega_2, \omega_3) = (2,2,0)$;\qquad
2)\ \ 
$|E(W)| =2$, $(\omega_1, \omega_2, \omega_3) = (2,0,2)$;\medskip\\
3)\ \ 
$|E(W)| =2$, $(\omega_1, \omega_2, \omega_3) = (2,1,1)$;\qquad
4)\ \ 
$|E(W)| =2$, $(\omega_1, \omega_2, \omega_3) = (1,2,1)$;\medskip\\
5)\ \ 
$|E(W)| =1$, $(\omega_1, \omega_2, \omega_3) = (2,2,1)$;\qquad
6)\ \ 
$|E(W)| =1$, $(\omega_1, \omega_2, \omega_3) = (2,1,2)$;\medskip\\
7)\ \ 
$|E(W)| =0$, $(\omega_1, \omega_2, \omega_3) = (2,2,2)$.
\end{quote}

In the case 1) or 5),
it follows from Lemma \ref{5cycle-lm} that
$\Xi^*-\{w_3\}$ is collapsible.
If $d_H(w_3) \ge 2$,
then we can find a DCT of $H$
by applying Lemma \ref{nokori2hen} to the graph $H/(\Xi^*-\{w_3\})$ and the vertices $v_{\Xi^*-\{w_3\}}$ and $w_{3}$.
If $d_H(w_3) = 1$,
then since the neighbor of $w_3$ is in $\Xi^*-\{w_3\}$,
the spanning closed trail of $\Xi^*-\{w_3\}$ is a DCT of $H$.
In the case 2) or 3),
since $\Xi^*$ is triangle-free,
we can take $u_i \in N(w_1)$ and $u_j \in N(w_3)$ so that
$u_i \in \{u_j\} \cup N_C(u_j)$.
In the case $u_i = u_j$ (resp.~$u_i \in N_C(u_j)$),
let $C' = u_i w_1 w_2 w_3 u_i$ (resp.~$C' = u_i w_1 w_2 w_3 u_j u_i$),
then $C \triangle C'$ is a spanning closed trail of $\Xi^*$.
In the case 4), without loss of generality we may assume that
$N_C(w_1)= \{u_1\}$.
Then, by symmetry, we may assume that
$N_C(w_3)= \{u_i\}$ for $i=1,2,3$.
If $N_C(w_3)= \{u_1\}$ (resp.~$u_2$),
then $C \cup u_1 w_1 w_2 w_3 u_1$ 
(resp.~$(C-\{u_1 u_2\}) \cup u_1 w_1 w_2 w_3 u_2$)
is a spanning closed trail of $\Xi^*$.
If $N_C(w_3)= \{u_3\}$, then by symmetry
we may assume that $N_C(w_2)= \{u_2, u_4\}$,
and then $E(\Xi^*) - \{u_1 u_2, u_3 u_4\}$ induces a spanning closed trail of $\Xi^*$.

\begin{figure}[tb]
\begin{center}
\begin{picture}(120,100)
\put(45,100){\footnotesize$x_2$}
\put(143,100){\footnotesize$x_3$}
\put(-14,65){\footnotesize$w_1$}
\put(45,65){\footnotesize$w_2$}
\put(143,65){\footnotesize$w_3$}
\put(-8,24){\footnotesize$u_1$}
\put(27,24){\footnotesize$u_2$}
\put(61,24){\footnotesize$u_3$}
\put(95,24){\footnotesize$u_4$}
\put(143,24){\footnotesize$u_5$}
\includegraphics[width=5cm]{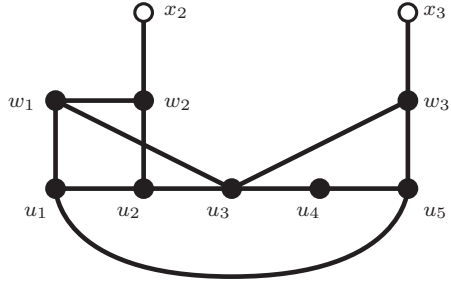}
\end{picture}
\caption{The graph $\Xi^*$ and two vertices $x_2$, $x_3$,
where black vertices denote $V(\Xi^*)$}
\label{fig_Xi}
\end{center}
\end{figure}

In the case 6),
without loss of generality we may assume that $u_2 \in N(w_2)$.
If $u_i \in N(w_1)$ for $i=4$ or $5$,
let $\Xi^{**}$ be the graph obtained from $\Xi^*-\{w_1, w_2\}$
by adding a vertex $w'$ and two edges $w' u_2, w'u_i$.
Then $\Xi^{**}$ is a triangle-free graph of order $7$ such that $V(\Xi^{**}) \setminus V(C) = \{w', w_{3}\}$, and hence
it follows from Lemma \ref{5cycle-lm} that
$\Xi^{**}$ contains a spanning closed trail $T$ such that
$\{w' u_2, w'u_i\} \in E(T)$.
Then $(T - w') \cup \{u_2 w_2 w_1 u_i\}$ is a spanning closed trail of $\Xi^*$.
Therefore, we may assume that $u_4, u_5 \notin N(w_1)$,
which implies $\{u_1, u_3\} \subseteq N(w_1)$ since $\Xi^*$ is triangle-free.
Then, by symmetry,
we may assume that
$\{u_1, u_3\}, \{u_2, u_4\}$ or $\{u_3, u_5\}\subseteq N(w_3)$.
If $\{u_1, u_3\}$ (resp.~$\{u_2, u_4\}$) $\subseteq N(w_3)$,
then $\Xi^* - \{w_1 u_3, u_2 u_3\}$ (resp.~$\Xi^* - \{w_1 u_1, u_3 u_4\}$) is a spanning closed trail of $\Xi^*$.
Hence we assume that $\{u_3, u_5\} \subseteq N(w_3)$
(See Figure \ref{fig_Xi}).

For $i=2,3$,
we may assume that $w_i$ has a neighbor $x_i$ in $H-\Xi^*$,
since $\Xi^* - \{w_i\}$ has a spanning closed trail $T_i$,
which is not a DCT of $H$.
If $x_2 u_5 \in E(G)$, then $u_5 x_2 w_2 w_1 u_3 u_2 u_1 u_5 u_4 u_3 w_3 u_5$
is a closed trail of $H$ which contains $V(\Xi^*)$,
and hence it is a DCT of $H$.
Thus we may assume that $x_2 u_5 \notin E(G)$.
Moreover, if there exists a path $P \subseteq H-\{u_1, u_3, u_4\}$
which joins $w_1$ and $w_3$ and contains $W$,
then $C \cup P \cup  \{u_3 w_1, u_3 w_3 \}$ is a DCT of $H$.
Hence there is no such a path, and thus
$x_{2} \neq x_{3}$ and
\begin{eqnarray}
E(\{x_2, w_2, u_2\}, \{x_3, w_3, u_5\}) = \emptyset.
\label{Nov30}
\end{eqnarray}

Now we have two subdivided claws induced by
$\{u_2 w_2 x_2, u_2 u_3 u_4, u_2 u_1 u_5\}$ and
$\{u_5 w_3 x_3,u_5 u_1 u_2, u_5 u_4 u_3\}$.
Let $y_2 = u_2$ and $y_3 = u_5$,
then by the fact that $n \ge 15$ and Lemma \ref{endmove} i) and ii)
\footnote{To be precise,
we have to consider the case where there exists $e_1 \in l_H(w_i)$ with $d^e_H(e_1) \ge 5$
or $e_2 \in l_H(y_i)$ with $d^e_H(e_2) \ge 7$.
Since the proof of these cases is exactly the same as the proof written in the above,
we omitted it.
In the rest of this paper,
we make the same omission in order to improve the readability.},
we can deduce that
either $d^e_H(x_i w_i) \ge 5$ or $d^e_H(w_i y_i) \ge 7$ holds
for $i=2$ and $3$.
Whichever inequality holds,
since $d_{\Xi^*}(w_i) = 2$ and $d_{\Xi^*}(y_i)=3$,
we can find $E_i \subseteq E(H) \setminus E(\Xi^*)$
such that $|E_i| \ge 4$ and each edge of $E_i$ has at least one endvertex in $\{x_i, w_i, y_i\}$
for $i=2$, $3$
(notice that $x_i w_i \in E_i$).
By (\ref{Nov30}), $E_2 \cap E_3 = \emptyset$,
and hence $n = |E(H)| \ge |E(\Xi^*)| + |E_2| + |E_3| \ge 19$,
which contradicts (\ref{Dec14}).

In the case 7),
since $\Xi^*$ is triangle-free,
it follows for each $w \in W$ that
$N(w) = \{u_l, u_{l+2}\}$
for some $l$ with $1 \le l \le 5$
(where indices are taken modulo $5$).
If $N(w) \neq N(w')$ for each $w, w' \in W$,
then we can take distinct $i,j \in \{1,2,3\}$ so that
$N(w_i) = \{u_l, u_{l+2}\}$ and $N(w_j) = \{u_{l+1}, u_{l+3}\}$
for some $l$.
If $N(w) = N(w') = \{u_l, u_{l+2}\}$ for some $w, w' \in W$ and $l$,
then since $\Xi^*$ has a perfect matching,
$u_{l+1} \in N(w'')$, where $w'' \in W \setminus \{w,w'\}$.
Consequently, we can take $w_i, w_j \in W$ so that
$N(w_i) = \{u_l, u_{l+2}\}$ and $u_{l+1} \in N(w_j)$ for some $l$.
Let $w_k \in W \setminus \{w_i, w_j\}$.
By Lemma \ref{5cycle-lm},
$\Xi^*-\{w_i\}$ contains a spanning closed trail $T$ such that $w_j u_{l+1} \in E(T)$.
Let $C'$ be the cycle $w_i u_l u_{l+1} u_{l+2} w_i$,
then $T \triangle C'$ is a spanning closed trail of $\Xi^*$.
\medskip\\
%
%
%
%
%
%
\noindent
\textbf{Case 2.}\quad
$n \le 14$.
\medskip

If $n \le 8$, then since $G$ is $2$-connected,
minimum degree of $G$ is at least $2 \ge \frac{n-2}3$.
Hence by Theorem \ref{MS}, $G$ is hamiltonian.
Therefore we may assume that $9 \le n \le 14$.
Take a closed trail $T \subseteq H$ with an arbitrary direction
so that
$T$ dominates as many edges as possible.
We may assume that there is no DCT in $H$,
and hence $E(H-V(T)) \neq \emptyset$.
Let $S$ be the component of $H-V(T)$ containing
at least one edge of $E(H-V(T))$ and
let $U = H - (V(T) \cup V(S))$.
Moreover, let $X=\{x_1, x_2, \ldots ,x_k\}$ be the set of the vertices of $T$
which have a neighbor in $S$,
where $x_1, x_2, \ldots $ appear in this order along $T$
(if some $x_i$ appears twice or more in $T$,
then we focus on the first appearance hereafter),
and let $x_{k+1} = x_1$.
Then we can take edge-disjoint trails $T_1, \ldots , T_k \subseteq T$,
where $T_i$ is a segment of $T$ 
from $x_i$ to $x_{i+1}$ along the direction of $T$. 
Let $f$ be the map from $E(U,T-X)$ to $\{1,2, \ldots ,k\}$
such that $f(e) = \min\{ i \mid \mbox{$e$ is incident to a vertex in $V(T_i)\setminus \{x_i, x_{i+1}\}$} \}$,
and let $Y_i = \{ e \mid e \in E(U,T-X),\ f(e) = i\}$.
Moreover, let $F_i = E(T_i) \cup Y_i$ 
and $X_i = E(U, x_i)$.
Then we have $F_i \cap F_j = X_i \cap X_j = \emptyset$ for $i \neq j$.
In addition, $F_i \cap X_j = \emptyset$ for $1 \le i, j \le k$.

By the maximality of $T$,
each $x_i$ has exactly one neighbor in $S$.
Since $H$ is essentially $2$-edge-connected,
we have $k \ge 2$.
For $1 \le i \le k$,
let $\mathcal{P}_i$ be the set of paths of $H$
joining $x_i$ and $x_{i+1}$ whose internal vertices are contained in $S$,
and take $P_i \in \mathcal{P}_i$ so that $|P_i|$ is maximum.
Let $S_i = E(P_i) \cup \{ e \in E(S) \mid \mbox{$e$ is dominated by $P_i$}\}$,
then we have $|S_i| \ge 3$
because $S_i$ contains two edges of $E(S,T)$ and at least one edge of $E(S)$.
Note that $S_i \cap F_j = S_i \cap X_j = \emptyset$ for $1 \le i, j \le k$.

It follows for each $i$ that $|F_i| \ge |S_i|$,
since otherwise the closed trail induced by $(E(T) \cup E(P_i)) \setminus E(T_i)$
dominates more edges than $T$.
Therefore,
\[
14 \ge n = |E(H)| \ge \sum_{1 \le i \le k}|F_i| + |E(S,T)| + |E(S)|
\ge 3k + k + 1,
\]
which yields $k \le 3$.
For each $i$,
let $x_i^+$ (resp.~$x_i^-$) be the successor (resp.~predecessor) of $x_i$ in $T$, and
$x_i^{+2}$ (resp.~$x_i^{-2}$) be the successor (resp.~predecessor) of $x_i^+$ (resp.~$x_i^-$) in $T$.
Moreover, let $x_i'$ be a neighbor of $x_i$ in $S$ and let $x_i''$ be a neighbor of $x_i'$ in $S$.
Note that, by the choice of $T$,
neither $x_i^+$ nor $x_i^-$ has a neighbor in $S$ for each $i$.

\begin{claim}
$\{u x_i^-, u x_{i+1}^-\} \not \subseteq N_H(u)$ for each $u$ and $i$ with $u \in U$ and $1 \le i \le k$.
\label{tasuki}
\end{claim}
\noindent
\textit{Proof.}\quad
Assume to the contrary that $u x_i^-, u x_{i+1}^- \in E(H)$ for some $u \in U$.
Then  the closed trail induced by
$( E(T) \cup E(P_i) \cup \{u x_i^- , u x_{i+1}^-\} ) \setminus \{x_i^- x_i, x_{i+1}^- x_{i+1}\}$
contradicts the choice of $T$.
\qed
\medskip

For each $i$,
if there exists $u \in U$ with $u x_i^- \in Y_{i-1}$
(where indices are taken modulo $k$),
then let
$y_i = u$,
and if $u x_i^- \notin Y_{i-1}$ for every $u \in U$,
then let $y_i = x_i^{-2}$.
Note that $|F_i|\ge 3$ yields $y_i \neq x_{i-1}$.
\medskip
\begin{claim}
$|E(T)| \ge 5$.
\label{notC4}
\end{claim}
\noindent
\textit{Proof.}\quad
Assume to the contrary that $|E(T)| = 4$,
then $T = x_1 x_1^+ x_2 x_2^+ x_1$ and
$y_i \in U$ for $i=1,2$.
Let $a_1 = x_1 x_1^-$, $a_2 = x_1 x_1^+$, $a_3 = x_1 x_1'$,
$b_1 = x_2 x_2^+$, $b_2 = x_2 x_2^-$ and $b_3 = x_2 x_2'$. 
Moreover, let
$a_1^+ = b_1^+ = y_1 x_1^-$, $a_2^+ = b_2^+ = y_2 x_2^-$,
$a_3^+ = x_1' x_1''$ and $b_3^+ = x_2' x_2''$.

If $N_H^e(a_i^+) \cap N_H^e(a_j^+) \neq \emptyset$ for some $i \neq j$,
then
we can find a closed trail which contains $V(T)$ and a vertex in $S$,
which contradicts the choice of $T$.
Hence
$N_H^e(a_i^+) \cap N_H^e(a_j^+) = \emptyset$ for each $i,j$ with $i \neq j$.
By symmetry, we have
$N_H^e(b_i^+) \cap N_H^e(b_j^+) = \emptyset$ for each $i,j$ with $i \neq j$.
Hence each of $\{a_1, a_1^+, a_2, a_2^+, a_3, a_3^+\}$ and
$\{b_1, b_1^+, b_2, b_2^+, b_3, b_3^+\}$ induces a subdivided claw.

By Lemma \ref{endmove} i) and ii),
we can deduce that 
$\ed(a_i) \ge \frac{n-2}3+2$ or $\ed(a_i^+) \ge \frac{n-2}3$ for each $i$. 
If all of $a_1^+, a_2^+$ and $a_3^+$ are heavy, 
then we have $|E(H)| \ge 3 \times \frac{n-2}3 + 3 >n$, a contradiction.
If $a_i^+$ is not heavy for some $i$ and $a_j^+$ is heavy for each $j \neq i$,
then since $|N_H^e(a_i) \cap N_H^e(a_j^+)| \le 1$
for each $j \neq i$, 
we have $|E(H)| \ge  (\frac{n-2}3+2) + 2 \cdot \frac{n-2}3  + 3  - 1 - 1>n$,
a contradiction.
Moreover, if $a_i^+$, $a_j^+$ is not heavy 
and $a_k^+$ is heavy for $\{i,j,k\} = \{1,2,3\}$,
then we may assume $i \neq 3$ without loss of generality.
By the fact $a_i^+ = b_i^+$, we obtain $\ed(b_i) + \ed(a_j) + \ed(a_k^+) \ge (\frac{n-2}3+2) + (\frac{n-2}3+2) + \frac{n-2}3  \ge n+2$, which contradicts Lemma \ref{jisuukeisan}. 
These arguments imply that $a_i^+$ is not heavy for each $i$,
and by symmetry it follows that $b_i^+$ is not heavy for each $i$.
Then we can apply Lemma \ref{endmove} with $|J|=2$,
and hence 
$\ed(a_1), \ed(b_2) \ge \frac{n-2}3+4$ holds. 
Since $|N_H^e(a_1) \cap N_H^e(b_2)| \le 2$, 
we have 
$n = |E(H)| \ge |N_H^e(a_1)| + |N_H^e(b_2)| - |N_H^e(a_1) \cap N_H^e(b_2)| + |\{a_1, b_2\}|
\ge 2(\frac{n-2}3+4) -2 + 2 = \frac{2n + 20}{3}$.
This implies $n \ge 20$, a contradiction. 
\qed

\begin{claim}
$x_j^{-2} \neq x_j^{+2}$ for some $j$.
\label{puramai2}
\end{claim}
\noindent
\textit{Proof.}\quad
Assume to the contrary that $x_i^{-2} = x_i^{+2}$ for each $i$,
then $C_i = x_i^{-2} x_i^- x_i x_i^+ x_i^{+2}$ induces a cycle of length $4$ for each $i$.

If $x_2^{-2} = x_1$, then by the choice of $T$,
both of $x_2^-, x_2^+$ have a neighbor in $U$
(consider the trail $(T \cup P_1) - \{x_2^{-2} x_2^-, x_2^- x_2\}$ for example). 
Recall that $x_1^{-2} = x_1^{+2}$ and $x_2^{-2} = x_2^{+2}$.
If $x_2 \neq x_1^{+2}$, then 
$T = x_1^{-2} x_1^- x_1 x_1^+ x_1^{+2} \cdots x_1 x_2^- x_2 x_2^+ x_1 \cdots x_1^{-2}$, 
which implies $|E(T)| \ge 12$ (note that $H$ is triangle-free).Then $n \ge |E(T)| + |E(S, T)| + |E(S)|\ge 12 + 2 + 1 = 15$, a contradiction.  
Hence $x_2 = x_1^{+2} = x_1^{-2}$. 
Since $|E(T)| \ge 5$ follows from Claim \ref{notC4},
we have $x_1^- \neq x_2^+$ and $x_1^-$ have a neighbor in $U$
(consider the trail $(T \cup P_1) - \{x_1^{-2} x_1^-, x_1^- x_1\}$).
Note that $T' = x_1^{-2} x_1^- x_1 x_1^+ x_2 x_2^+ x_1$ is a subtrail of $T$.
This yields $|E(T)| \ge 8$, and hence $n = |E(H)| \ge |E(T)| + |E(U,\{x_1^-, x_2^-, x_2^+ \})| + |E(S,T)| + |E(S)|
\ge 8+3+2+1=14$.
This implies $|E(T)| = 8$, and hence $|E(T) \setminus E(T')| = 2$.
Let $x_2^{+3}$ be the vertex of $T$ such that $T = T' \cup x_2^{+2} x_2^{+3} x_1^{-2}$.
Then by the calculation above, $x_2^{+3}$ has no neighbor in $U$,
and hence the trail $T' \cup P_1$ contradicts the choice of $T$.
Therefore $x_2^{-2} \neq x_1$,
By the same reason, we have $x_{i+1}^{-2} \neq x_i$ for each $i$.

If $x_2^{-2} = x_1^+$, then
$T_1 = x_1^{-2} x_1^- x_1 x_1^+ x_1^{-2} x_2 x_2^+ x_1^+$
is a subtrail of $T$.
Since $x_1^+$ and $x_1^{-2}$ are adjacent, $|E(T) \setminus E(T_1)| \ge 3$.
Moreover, by the choice of $T$,
each of $x_1^-$, $x_2^+$ has a neighbor in $U$
(consider the trail $(T\cup P_1)-\{x_1^{-2} x_2, x_1^{-2} x_1^-, x_1^- x_1\}$ for example).
Then $n = |E(H)| \ge |E(T)| + |E(U,\{x_1^-, x_2^+ \})| + |E(S,T)| + |E(S)|
\ge (7+3)+2+2+1=15$, a contradiction.
This implies that $x_{i+1}^{-2} \neq x_i^{+}$ for each $i$,
and then we can deduce that $C_1, \ldots , C_k$ are edge-disjoint.

If $k = 3$, then we can deduce that
$n \ge |E(T)| + |E(S,T)| + |E(S)| \ge 4 \times 3 + 3 + 1 \ge 16$,
a contradiction.
If $k=2$ and $x_2^{-2} \neq x_1^{+2}$,
then $T-E(C_1 \cup C_2)$ contains two paths joining $x_2^{-2}$ and $x_1^{+2}$.
Since $H$ is triangle-free, the sum of the number of edges of these paths is at least $4$,
and $n \ge |E(T)| + |E(S,T)| + |E(S)| \ge (4 \times 2 + 4) + 2 + 1 \ge 15$,
a contradiction. Thus we have $k=2$ and $x_2^{-2} = x_1^{+2}$.

By the choice of $T$,
either $x_1^+$ or $x_2^-$ has a neighbor in $U$,
since otherwise  $(T \cup P_1) - \{x_1^+, x_2^-\}$ dominates an edge of $S$ and all the edges of $T$.
Similary, either $x_1^-$ or $x_2^+$ has a neighbor in $U$,
either $x_1^+$ or $x_2^+$ has a neighbor in $U$,
and
either $x_1^-$ or $x_2^-$ has a neighbor in $U$.
Therefore, without loss of generality we may assume that
each of $x_1^-$, $x_1^+$ has a neighbor in $U$.
Then $n \ge |E(T)| + |E(U,T)| + |E(S,T)| + |E(S)| \ge 8 + 2 + 2 + 1 \ge 13$.
Let $y^-$ be the neighbor of $x_1^-$ in $U$.
Then $\{x_1 x_1', x_1' x_1'', x_1 x_1^+, x_1^+ x_1^{+2}, x_1 x_1^-, x_1^- y^- \}$ induces a subdivided claw.
By Lemma \ref{endmove} i) and ii),
$\ed(x_1^- y^-) \ge 4$ or $\ed(x_1 x_1^-) \ge 6$ follows.
In either case, we have $|E(U,T)| + |X_1| \ge 4$.
This yields $n \ge 15$, a contradiction.
\qed
\begin{claim}
$k = 2$.
\end{claim}
\noindent
\textit{Proof.}\quad
Assume to the contrary that $k=3$.
Let $F_i' = \{x_i x_i^+, x_{i+1} x_{i+1}^-, x_{i+1}^- y_{i+1}\}$
(note that $F_i' \subseteq F_i$).
Moreover, let $E_1 = F_1' \cup F_2' \cup F_3' \cup E(S,T) \cup E(S)$,
then $|E_1| \ge 3 \cdot 3 + 3 + 1 = 13$.
Take $j$ as in Claim \ref{puramai2},
then since $H$ is triangle-free,
$\{x_j x_j', x_j' x_j'', x_j x_j^+, x_j^+ x_j^{+2}, x_j x_j^-, x_j^- y_j\}$ induces a subdivided claw.
Since we already have $n \ge |E_1| \ge 13$,
by Lemma \ref{endmove} i) and ii),
either $d^e(y_j x_j^-) \ge 4$ or $d^e(x_j^- x_j) \ge 6$ holds.
In either case, it follows from Claim \ref{tasuki} that
there exist two edges
in $E(H) \setminus E_1$
(note that, for $i \neq j$,
$x_j^- \neq x_i$ follows by the maximality of $T$, 
and if $x_j^- \in V(T_i)$ for $i \neq j-1$, then two closed trails
induced by $(E(T) \cup E(P_i)) \setminus E(T_i)$
and $(E(T) \cup E(P_{j-1})) \setminus E(T_{j-1})$
yield $|F_i|, |F_{j-1}| \ge 4$.) 
Hence $n \ge 13 + 2 = 15$, a contradiction.
\qed
\medskip

Assume that $x_i$ appears twice or more in $T_i$ for some $i$.
Let $T_i'$ be the segment of $T_i$ between the last appearance of $x_i$ and $x_{i+1}$,
then by the choice of $T$, $|E(T_i')| + |E(T_i'-\{x_i, x_{i+1}\}, U)| \ge |S_i| \ge 3$.
Let $E_2 = F_{i+1} \cup E(T_i) \cup E(T_i'-\{x_1, x_{i+1}\}, U) \cup E(S,T) \cup E(S)$,
then
since $E(T_i) \setminus E(T_i')$ contains a cycle of length at least $4$,
we have
\begin{eqnarray}
|E_2| \ge 3 + 4 + 3 + 2 + 1 = 13.
\label{eq_E2}
\end{eqnarray}
Since $H$ is triangle-free,
$\Lambda = \{x_i x_i', x_i' x_i'', x_i x_i^+, x_i^+ x_i^{+2}, x_i x_i^-, x_i^- y_i \}$
induces a subdivided claw.
If $d^e(x_i' x_i'') \ge 4$, $d^e(x_i^+ x_i^{+2}) \ge 4$ or $d^e(x_i^- y_i) \ge 4$,
then we have two edges which is not counted
in the inequality (\ref{eq_E2}),
even when $x_i^+$, $x_i^{+2}$, $x_i^-$ or $y_i$ appears twice or more in $T$
Such edges yield $n \ge 15$, a contradiction.
Since we already have $n \ge |E_2| \ge 13$,
by Lemma \ref{endmove} i) and ii),
we have $d^e(x_i x_i') \ge 8$ holds
(note that we can apply Lemma \ref{endmove} with $|J| = 2$).
Again, we have two edges which is not counted in the inequality (\ref{eq_E2}),
which yields $n \ge 15$, a contradiction.

Therefore, each $x_i$ appers only once in $T$.
We distinguish two cases.
\medskip\\
\textbf{Case 2.1.}\quad
$k = 2$ and $|E(S)| \ge 2$.
\medskip\\
Note that $|S_i| \ge |E(S)|+2$ holds for $i=1,2$,
since otherwise an edge of $S$ is separated by
another edge of $S$,
which contradicts the assumption that
$H$ is essentially $2$-edge-connected.
Hence $14 \ge |E(H)| \ge |F_1|+|F_2|+|E(S,T)| + |E(S)|
\ge |S_1|+|S_2|+2 + |E(S)| \ge 3|E(S)|+6$,
which implies $|E(S)| \le 2$.

Take $j$ as in Claim \ref{puramai2}.
Since $H$ is triangle-free,
$\Lambda = \{x_j x_j', x_j' x_j'', x_j x_j^+, x_j^+ x_j^{+2}, x_j x_j^-, x_j^- y_j\}$ induces a subdivided claw.
Since we already have $|E(H)| \ge 3|E(S)|+6 \ge 12$,
by Lemma \ref{endmove} i) and ii),
either $d^e(x_j' x_j'') \ge 4$ or $d^e(x_j x_j') \ge 6$ holds.
Since $d^e(x_j' x_j'') \le |E(S) \setminus \{x_j' x_j''\}| + |E(S,T)| \le 3$,
we have $d^e(x_j x_j') \ge 6$.
Then $N_H^e(x_j x_j') \setminus (E(S) \cup \{x_j' x_{j+1}, x_j x_j^-, x_j x_j^+ \}) \neq \emptyset$,
which implies $|X_j| \ge 1$.

If $x_{j+1}^{-2} \neq x_{j+1}^{+2}$,
then by the same argument as above,
we obtain $|X_{j+1}| \ge 1$,
which implies $|F_1| + |F_2|+|X_{j+1}| \ge |S_1| + |S_2| +|X_{j+1}|\ge 4 + 4+1 = 9$.
On the other hand,
if $x_{j+1}^{-2} = x_{j+1}^{+2}$,
then by Claim \ref{notC4},
$T$ induces at least two cycles of length at least $4$.
Moreover, one of the two cycles which contains $x_j$ has length at least $5$,
since $x_j^{-2} \neq x_j^{+2}$.
Hence $|F_1| + |F_2| \ge |E(T)| \ge 9$.
Consequently, we have $|F_1| + |F_2| + |X_{j+1}| \ge 9$ in either case.
Let $E_3 = F_1 \cup F_2 \cup X_j \cup X_{j+1} \cup E(S,T) \cup E(S)$,
then $|E(H)| \ge |E_3| \ge 9 + |X_j| + |E(S,T)| + |E(S)| \ge 14$,
which implies $|F_1| = |F_2| = 4$ and $|X_j| = |X_{j+1}|=1$ in the case $x_{j+1}^{-2} \neq x_{j+1}^{+2}$,
and $|F_1|+|F_2| = |E(T)| = 9$, $|X_j| = 1$ and $|X_{j+1}| = 0$ in the case $x_{j+1}^{-2} = x_{j+1}^{+2}$.
In the former case, we obtain $d^e(y_j x_j^-) \le 3$,
and in the latter case,
we obtain $d^e(y_j x_j^-) \le 3$ (in the case $y_{j} \neq x_{j+1}^{-2}$)
or $d^e(x_j^{+2} x_j^+) \le 3$ (in the case $y_{j} = x_{j+1}^{-2})$).
Recall that $\Lambda$ induces a subdivided claw.
Since $d^e(y_j x_j^-)$ or $d^e(x_j^{+2} x_j^+)$ is at most $3$,
it follows from Lemma \ref{endmove} i) and ii) that $d^e(x_j x_j') \ge 7$
(note that we can apply Lemma \ref{endmove} with $|J| \ge 1$).
This yields $|X_j| \ge 2$, a contradiction.
\medskip\\
\textbf{Case 2.2.}\quad
$k = 2$ and $|E(S)| = 1$.
\medskip\\
Let $\alpha = 0$ in the case $9 \le n \le 11$,
and let $\alpha = 1$ in the case $12 \le n \le 14$.
Take $j$ as in Claim \ref{puramai2}.
Since $H$ is triangle-free,
$\Lambda = \{x_j x_j', x_j' x_j'', x_j x_j^+, x_j^+ x_j^{+2}, x_j x_j^-, x_j^- y_j\}$ induces a subdivided claw.
By Lemma \ref{endmove} i) and ii),
either $d^e(x_j' x_j'') \ge 3 + \alpha$ or $d^e(x_j x_j') \ge 5 + \alpha$ holds.
Since $d^e(x_j' x_j'') \le |E(S,T)| \le 2$,
we have $d^e(x_j x_j') \ge 5 + \alpha$.
Then $N_H^e(x_j x_j') \setminus \{x_j' x_j'', x_j' x_{j+1}, x_j x_j^-, x_j x_j^+ \} \neq \emptyset$,
which implies $|X_j| \ge 1 + \alpha$.

Let $E_4 = F_1 \cup F_2 \cup X_j \cup X_{j+1} \cup E(S,T) \cup E(S)$.
Assume first that $x_{j+1}^{-2} = x_{j+1}^{+2}$,
then by the same argument as in Case 2.2,
$|F_1| + |F_2| \ge |E(T)| \ge 9$.
Hence we have $|E_4| = 9 + 1 + \alpha + 2 + 1 = 13 + \alpha$.
Since $n \ge |E_4| \ge 13$, we have $\alpha = 1$,
and then the fact $|E(H)| \le 14$ yields
$|F_1| + |F_2| = |E(T)| = 9$,
$|X_j| = 2$ and
$|X_{j+1}| = 0$.
Since $|F_1| + |F_2| = |E(T)|$,
we have $y_j = x_j^{-2}$ and $d^e(y_j x_j^-) = 2$.
Since $\Lambda$ induces a subdivided claw,
it follows from Lemma \ref{endmove} i) and ii) that $d^e(x_j x_j') \ge 7$
(note that we can apply Lemma \ref{endmove} with $|J| \ge 1$).
This implies $|X_j| \ge 3$, a contradiction.

Assume next that $x_{j+1}^{-2} \neq x_{j+1}^{+2}$.
Then by the same argument as above,
we obtain $|X_{j+1}| \ge 1 + \alpha$.
Recall that $\Lambda$ induces a subdivided claw.
If $d^e(y_j x_j^-) \ge 3 + \alpha$,
then we can deduce that $|F_1| + |F_2| + |E(H) \setminus E_4| \ge 7+\alpha$
(note that, if an edge $e$ is adjacent to $y_j x_j^-$,
then either $e \in F_1 \cup F_2$ or $e \in E(H) \setminus E_4$).
Hence $|E(H)| \ge |F_1| + |F_2| + |E(H) \setminus E_4| + |X_1| + |X_2| + |E(S,T)| + |E(S)|
 \ge  7+\alpha + 2(1+\alpha) + 2+1 = 12 + 3\alpha$,
which is a contradiction in either case of $n \le 11$ or $n \ge 12$.
Hence $d^e(y_j x_j^-) < 3 + \alpha$, which implies
$d^e(x_j x_j') \ge 6 + \alpha$
(note that we can apply Lemma \ref{endmove} with $|J| \ge 1$).
Thus we have $|X_j| \ge 2 + \alpha$.
Since $x_{j+1}^{-2} \neq x_{j+1}^{+2}$,
$\Lambda = \{x_{j+1} x_{j+1}', x_{j+1}' x_{j+1}'', x_{j+1} x_{j+1}^+, x_{j+1}^+ x_{j+1}^{+2}, x_{j+1} x_{j+1}^-, x_{j+1}^- y_{j+1}\}$
induces a subdivided claw.
Hence, by the same argument as above,
we obtain $|X_{j+1}| \ge 2 + \alpha$.
Thus $|E(H)| \ge |F_1| + |F_2| + |X_1| + |X_2| + |E(S,T)| + |E(S)|
\ge 3 + 3 + 2(2+\alpha) + 2 + 1 = 13 + 2\alpha$,
which is a contradiction in either case of $n \le 11$ or $n \ge 12$.
\qed
%
%
%
%
%
%
%

%
%
%
\end{document}